\input amstex
\define \g {\gamma}
\define \db {\bar{\partial}}
\define \Ker {\operatorname{Ker}}
\define \C {\Bbb C}
\define \Def {\text {Def}}
\define \p {\partial}
\define \gt {\goth}
\define \gtg {\goth g}
\define \Z {\Bbb Z}
\define \om {\Omega}
\define \dv {\Delta}
\define \dl {\Delta}
\define \im {\operatorname{Im}}
\define \Hom {\operatorname{Hom}}
\define \ta {\tau}
\define \bt {\bold t}

\define\proof{\demo{Proof}}
\define\endproof{\qed\enddemo}
\define\theorem#1.#2{\proclaim{Theorem #1.#2}}
\define\lemma#1.#2{\proclaim{Lemma #1.#2}}
\define\proposition#1.#2{\proclaim{Proposition #1.#2}}
\define\corollary#1.#2{\proclaim{Corollary #1.#2}}
\define\claim#1.#2{\proclaim{Claim #1.#2}}
\define\Remark{\remark{Remark}}

\define\section#1{\specialhead #1 \endspecialhead}
\define\subsection#1.#2{\subhead #1.#2 \endsubhead}

\documentstyle{amsppt}

\topmatter
\title Frobenius manifolds and Formality of Lie algebras of polyvector fields
\endtitle
\rightheadtext{Frobenius manifolds and Formality of Lie algebras}
\author {Sergey Barannikov, Maxim Kontsevich}
\endauthor
\address University of  California at Berkeley, Berkeley CA 94720, 
USA\endaddress
\address Institut des Hautes \'Etudes Scientifiques, 
35 route de Chartres, 91440
Bures-sur-Yvette, France
\endaddress
\abstract
We construct a generalization of the variations of Hodge structures
on Calabi-Yau manifolds. It gives a Mirror
 partner for the theory of $genus=0$ 
Gromov-Witten invariants. 
\endabstract
\email barannik\@ihes.fr , maxim\@ihes.fr\endemail
\thanks  S.~B.~ was partially supported by the Fellowship for
Graduate Study of the University of California at Berkeley\endthanks

\endtopmatter

\vskip 2truecm

\nologo

\document

\section { Introduction}

Probably the best mathematically understood part of the 
Mirror Symmetry program is the theory of Gromov - Witten invariants (see [KM]).
In this paper we will construct a Mirror partner for the $genus=0$ sector of
 this theory. 
It may be 
considered as a generalization  of the theory of variations of Hodge structures
on Calabi-Yau manifolds.

One of the puzzles in Mirror symmetry  was 
to find an interpretation of the mysterious objects 
involved in the famous predictions of the numbers of
rational curves.  
Such an interpretation should, in particular, give the meaning to the
``extended'' moduli space  $H^*(M,\Lambda^*T_M)[2]$, thought as 
generalized deformations
of complex structure. This moduli space should be equipped with the analog
of the
3-tensor $C_{ijk}(t)$ (``Yukawa coupling'') arising from a
generalization of the variation of 
Hodge structure on $H^*(M)$.  To find such structure is essential for
the extension of the predictions of Mirror Symmetry in the dimensions $n>3$.

\subhead 0.1 Background philosophy\endsubhead

The Mirror Symmetry conjecture, as it was formulated in [K1], states that the
derived category of coherent sheaves on a Calabi-Yau manifold $M$ is 
equivalent
to the derived category constructed from (conjectured) Fukaya category 
associated 
with the
dual Calabi-Yau manifold $\widetilde M$.  The conjecture implies existence of
the structure of Frobenius manifold on the moduli space of
$A_{\infty}$-deformations of the derived category of coherent sheaves on $M$.
This structure coincides conjecturally with the Frobenius structure on formal
neighborhood of zero in $H_*(\widetilde M,\C)$ constructed via
Gromov-Witten classes of the dual Calaby-Yau manifold $\widetilde M$.

\subhead 0.2 Formulation of the results\endsubhead

Consider the differential graded Lie algebra $$\bt=\bigoplus_k \bt^k,
\,\bt^k=\bigoplus_{q+p-1=k} 
\Gamma(M,\Lambda^q\overline{T}^*_M\otimes\Lambda^pT_M)
\eqno (0.1)$$ endowed with the differential $\db$ and the bracket 
coming from the cup-product on $\db$-forms and standard
Schouten-Nijenhuys bracket on holomorphic polyvector fields.
Deformation theory associates a formal (super)moduli space $\Cal M_{\bt}$ 
to the
Lie algebra $\bt$.  It can be described informally as universal moduli space of
solutions to the Maurer-Cartan equation over $\Z$-graded Artin algebras modulo
gauge equivalence.  The formal moduli space $\Cal M_{\bt}$ can be identified 
with
a formal neighborhood of zero in graded vector space $H^*(M,\Lambda^*T_M)[2]$ 
\footnote {For a graded object $\bold t$ denote
 $\bold t[n]$ 
the tensor product of $\bold t$ with the trivial object concentrated
in degree $-n$.}. 
 The tangent sheaf
to $\Cal M_{\bt}$ after the shift by $[-2]$ has natural structure of graded 
commutative 
associative algebra 
over $\Cal O_{\Cal M_\bt}$.
In this note we show that this algebra structure gives rise to the structure
 of
(formal) Frobenius manifold on $H^*(M,\Lambda^*T_M)[2]$.  More specifically,
using a universal solution to the analog of Maurer-Cartan equation in 
$\bt$, we
construct generalized holomorphic volume element for generalized
 deformations of
complex structure.  The integrals of this element, which can be thought of as
generalized periods, produce the Frobenius manifold structure on
$H^*(M,\Lambda^*T_M)[2]$.

\subhead 0.3 Connection with  $A_{\infty}$-deformations of 
$\Cal D^bCoh(M)$\endsubhead 

The Formality theorem (see [K2])  identifies the germ of the 
 moduli space of $A_{\infty}$-deformations of the
derived category of coherent sheaves on $M$ with the moduli space
 $\Cal M_{\bt}$. 
The tangent bundle of this moduli space after the shift by $[-2]$ 
has natural structure of the graded 
commutative associative algebra. The multiplication arises from the 
Yoneda product on $\text{Ext}$-groups. The identification of moduli spaces 
provided by the Formality theorem respects the algebra structure on
the tangent 
bundles of the moduli spaces.  This implies, in particular, that the 
usual
predictions of numerical Mirror Symmetry can be deduced  from the homological 
Mirror Symmetry conjecture  proposed in [K1]. We hope to elaborate on this
elsewhere.

\section{1. Frobenius manifolds}

 Remind the   definition of formal Frobenius (super) manifold as given in [D], 
[M], 
[KM].  Let $\bold H$ be a finite-dimensional $\bold Z_2$-graded 
vector space over
 $\C$.\footnote{One can use an
 arbitrary field of characteristic zero instead of $\C$ everywhere} It 
is convenient to choose 
some set of coordinates $x_{\bold H}=\{x^a\}$ which defines 
the basis $\{\p_a:=\p/\p x^a\}$
of vector fields. One of the given coordinates is distinguished and 
is denoted by $x_0$.
 Let $A^c_{ab}\in \C [[x_{\bold H}]]$ be a formal power series 
 representing 3-tensor field,
$g_{ab}$ be a  nondegenerate  symmetric pairing on $\bold H$. 
To simplify notations in superscripts we replace 
$\text{deg}\, (x^a)$ by $\bar a$.

One can use the $A_{ab}^c$ in order to define a structure of 
$\C [[x_{\bold H}]]$-algebra
on $\bold H\otimes\C [[x_{\bold H}]]$, the (super)space of all
continuous derivations of 
$\C[[x_{\bold H}]]$,
 with multiplication denoted by $\circ$:
$$ \p_a\circ\p_b:=\sum_c A_{ab}^c\p_c$$
One can use $g_{ab}$ to define the symmetric $\C[[x_{\bold H}]]$-pairing
on $\bold H\otimes\C[[x_{\bold H}]]$:
$$\langle \p_a,\p_b\rangle:=g_{ab}$$

These data define the structure
of formal Frobenius manifold on $\bold H$ iff the following equations hold:
\item{(1)} (Commutativity/Associativity)
$$\forall a,b,c \,\,\,\,\,\, A^c_{ba}=(-1)^{\bar a\bar b}A_{ab}^c\eqno (1a)$$ 
$$\forall a,b,c,d\,\,\,\,\,\,\sum_{e}A^e_{ab} A_{ec}^d 
=(-1)^{\bar a(\bar b+\bar c)}\sum_{e}A_{bc}^e 
A_{ea}^d\eqno (1b)$$
equivalently, $A^c_{ab}$ are  structure constants of a supercommutative
associative \break $\C[[x_{\bold H}]]$-algebra  
\item {(2)}(Invariance) Put $A_{abc}=\sum_e A^e_{ab}g_{ec}$  $$\forall a,b,c 
\,\,\,\,\,\,A_{abc}= 
(-1)^{\bar a(\bar b+\bar c)}A_{bca}\,\,\,\, ,$$ equivalently,
the pairing $g_{ab}$ is invariant with respect to the multiplication $\circ$ 
defined by 
$A_{ab}^c$.

\item {(3)} (Identity) $$\forall a,b \,\,\,\,\,\,A_{0a}^b=\delta^b_a$$ 
equivalently 
$\p_0$ is an identity of 
the algebra $\bold H\otimes\C[[x_{\bold H}]]$
\item {(4)} (Potential) $$\forall a,b,c,d\,\,\,\,\,\, \p_d 
A^c_{ab}=(-1)^{\bar a\bar d}\p_a 
A^c_{db}\,\,\,\, ,$$  which implies, assuming (1a) and (2), that the 
series $A_{abc}$ are the third derivatives of a single power series
 $\Phi\in \bold H\otimes
\C[[x_{\bold H}]]$
$$A_{abc}=\p_a\p_b\p_c\Phi$$

\section{2. Moduli space via deformation functor}

The material presented in this section is standard (see [K2] and references
 therein).

Let us remind the definition of  the
functor ${\Def_\gtg}$
 associated with a differential graded Lie algebra $\goth g$ . 
It acts from the category of   Artin algebras  to the category of 
sets. 
Let $\goth A$ be an Artin algebra with the  maximal ideal denoted by $\gt m$. 
Define the set   $${\Def_\gtg}(\gt A):=\{d\g+{[\g,\g]\over 2}=0 |\g
\in (\gt g\otimes\gt m)^1\}/\Gamma_{\gt A}^0$$
where the quotient is taken modulo action
of the group  $\Gamma_{\gt A}^0$ corresponding to the nilpotent Lie algebra  
$(\gt g\otimes\gt m)^0$. The action of the group can be described via the 
infinitesimal action of its
Lie algebra: $$\alpha\in \gt g\otimes\gt m \to \dot \g=d\alpha +[\g,\alpha]$$

Sometimes functor  $\Def_\gtg$ is represented by
some topological algebra
$\Cal O_{\Cal M_\gtg}$ (projective limit of Artin algebras) in the sense that 
the functor 
$\Def_\gtg$ is equivalent to the functor 
$\text{Hom}_{continuous}(\Cal O_{\Cal M_\gtg},\,\cdot\,)$.
For example, $H^0(\gtg)=0$ is  a sufficient condition for  this.
If $\Def_\gtg$ is representable then one can associate  formal moduli space to 
$\gtg$ by defining   
the  ``algebra of functions''
on the formal  moduli space to be  the \break algebra $\Cal O_{\Cal M_\gtg}$.

We will need the $\Z$-graded extension of the functor $\Def_\gtg$.
  The definition
of  $\Def^{\,\Z}_\gtg$ is obtained from the definition of 
$\Def_\gtg$ via inserting
$\Z-$graded Artin algebras  instead of the usual ones everywhere.
A sufficient and probably necessary condition for the functor 
$\Def^{\,\Z}_\gtg$ to be representable is that $\gtg$ must be 
 quasi-isomorphic to an abelian graded Lie algebra.
We will see in \S 2.1 that this is the case for $\gtg=\bt$.
Hence one can associate formal (graded) moduli space \footnote{We will
omit the superscript $\Z$ where it does  not seem to lead to a confusion.}
$\Cal M_\bt$ 
to the Lie algebra $\gtg$.

\subhead {2.1 Extended moduli space of   complex structure}
\endsubhead

Let $M$ be a connected compact complex manifold of dimension $n$, 
with vanishing 
$1$-st
Chern class $c_1(T_M)=0\in \text{Pic}(M)$. We assume that there exists a 
K\"ahler metric on $M$, although we will not fix it. By Yau's theorem there 
exists a Calabi-Yau metric on $M$.

It follows from the condition $c_1(T_M)=0$ there exists an everywhere 
nonvanishing holomorphic volume form $\Omega\in \Gamma(X,\Lambda^n T_M^*)$.
It is defined up to a multiplication by a constant. Let us fix a choice of 
$\om$.
It induces isomorphism
of cohomology groups: $H^q(M,\Lambda^pT_M)\simeq H^q(M,\Omega^{n-p})$;
$\g\mapsto \gamma\vdash\om$.
Define differential $\dv$ of degree $-1$ on $\bt$ by the formula :
$$(\dv\g)\vdash\om=\p(\g\vdash\om)$$ 
The operator $\dv$ satisfies the following identity (Tian-Todorov lemma) :
$$ [\g_1,\g_2]=(-1)^{\text{deg}\g_1+1}(\dv(\g_1\wedge\g_2) -
(\dv \g_1)\wedge\g_2 -
(-1)^{\text{deg}\g_1+1}\g_1\wedge\dv\g_2)\eqno (2.1)$$
where $\text{deg}\g=p+q-1$ for $\g\in\Gamma(M,\Lambda^q\bar 
T^*_M\otimes\Lambda^q
T_M)$.

Denote by $\bold H$ the total homology space of
$\dv$ acting on $\bt[1]$. Let $\{\dl_a\}$ denote a graded basis in the
vector space $\oplus_{p,q}H^q(M,\Lambda^pT_M)$, 
$\dl_0=1\in H^0(M,\Lambda^0T_M)$  .
 Let us redefine the degree 
of $\dl_a$ as follows 
$$|\dl_a|:=p+q-2\,\,\,\,\,\text{for}\, \,\,\,\,\,\dl_a\in 
H^q(M,\Lambda^pT_M)$$
Then $\{\dl_a\}$ form a graded basis in $\bold H$. 
Denote
by $\{t^a\}, t^a\in \bold H^*,\,\text{deg}\,t^a=-|\dl_a|$  the  basis  
dual to $\{\dl_a\}$. Denote by 
$\C[[t_{\bold H}]]$ the algebra of formal power series  on 
$\Z$-graded vector space $\bold H$. 

\proclaim{Lemma 2.1} The functor $\Def^{\,\Z}_\bt$ associated with 
$\bt$ is canonically
 equivalent to the functor represented by  the
algebra  $\C[[t_{\bold H}]]$.
\endproclaim
\proof
  It follows from (2.1) that the maps  $$(\bt,\db)\gets(\Ker\dv,\db)\to (\bold 
H[-1],d=0)\eqno (2.2)$$ are morphisms of differential graded Lie 
algebras. Then the
$\p\db$-lemma, which says that
$$\Ker\db\cap \Ker\dv\cap(\im\dv\oplus\im\db)=\im\dv\circ\db,\eqno (2.3) $$
 shows that  these  
morphisms are
quasi-isomorphisms (this argument is standard in the theory of minimal models,
see [DGMS]). 
Hence (see e.g. theorem in \S 4.4 of [K2])  the deformation 
functors associated with the three differential graded Lie algebras are 
canonically equivalent. The deformation functor associated with trivial algebra 
$(\bold H[-1],d=0)$ is represented by the algebra  $\C[[t_{\bold H}]]$.
\endproof

\corollary 2.2 The moduli space $\Cal M_\bt$ associated to 
$\bt$ is smooth. The
 dimension  of $\Cal M_\bt$ is equal to  
$\sum_{p,q}\text{dim }H^q(M,\Lambda^pT_M)$ of the dimension of the space 
of first order deformations associated with $\bt$.\endproclaim
\Remark  The Formality theorem proven in [K2] implies 
that the differential graded Lie algebra controlling the 
$A_\infty$-deformations
of $\Cal D^bCoh(M)$ is  quasi-isomorphic to $\bt$. Here we have proved that 
$\bt$ is quasi-isomorphic to an abelian graded Lie algebra. Therefore, 
 the two differential graded Lie algebras are formal, 
 i.e. quasi-isomorphic to their cohomology Lie algebras
endowed with zero differential.\endremark

\corollary 2.3  There exists a  solution to the Maurer-Cartan equation 
$$\db{\hat\g}(t)+{[{\hat\g}(t),{\hat\g}(t)]\over 
2}=0\eqno (2.4)$$
in formal power series with values in $\bt$ 
$${\hat\g}(t)=\sum_a\hat\g_at^a+{1\over 
{2!}}\sum_{a_1,a_2}\hat\g_{a_1a_2}t^{a_1}t^{a_2} +\ldots\in (\bt\otimes\bold 
C[[t_{\bold H}]])^1$$ such that
the cohomology classes $[\hat\g_a]$ form a basis of cohomology of the 
complex 
$(\bt,\db)$\endproclaim

\Remark The  deformations of the complex structure are controlled by  the
differential graded Lie algebra 
$$\bt_{(0)}:=\bigoplus_k \bt_{(0)}^k, \,\,\bt_{(0)}^k=
\Gamma(M,\Lambda^{k}\overline{T}^*_M\otimes T_M)$$ The
 meaning of this is that
the completion of the algebra of functions on the moduli space of complex 
structures on $M$ represents $\Def_{\bt_{(0)}}$ (or $\Def_{\bt_{(0)}}^{\,\Z}$  
 restricted to the category of Artin algebras concentrated in degree 
$0$). The natural embeddings $\bt_{(0)}\hookrightarrow \bt$ induces
embedding of the corresponding moduli spaces. In terms of the solutions
to Maurer-Cartan equation the deformations of complex 
structure are singled out 
by the condition $\g(t)\in 
\Gamma(M,\Lambda^1\overline{T}^*_M\otimes\Lambda^1T_M)$.\endremark
\Remark Similar thickening of the moduli space of complex structures were 
considered by Z\.~Ran in [R].\endremark

\section{3. Algebra structure on the tangent sheaf of the moduli space }
Let   $R$ denotes a  $\Z$-graded Artin algebra over $\C$, $\g^R
\in (\bt\otimes R)^1$ denotes
a solution to the Maurer-Cartan equation (2.4). 

 The linear extension of the wedge product gives a  structure of graded 
commutative algebra on $\bt\otimes R[-1]$. Let $\g^R$ be a solution to the 
Maurer-Cartan equation in $(\bt\otimes R)^1$. It defines a differential 
$\db_{\g^R}=\db+\{\g^R,\cdot\}$ on $\bt\otimes R[1]$.  Denote the cohomology of 
$\db_{\g^R}$ by $T_{\g^R}$.
The space of first order variations of $\g^R$ modulo gauge equivalence  is 
identified with $T_{\g^R}$. Geometrically one can think of $\g^R$ as
 a morphism 
from the formal variety corresponding to algebra $R$ to the formal moduli 
space. An element of $T_{\g^R}$ corresponds to a section of the 
preimage of the 
tangent sheaf. 
Note that 
$\db_{\g^R}$ acts as a differentiation of the (super)commutative $R-$algebra 
$t\otimes R[-1]$.
Therefore  $T_{\g^R}[-2]$ inherits the structure of (super)commutative 
associative
algebra over $R$. This structure is functorial with respect to the morphisms
$\phi_*: T_{\g^{R_1}}\to T_{\g^{R_2}}$ induced by homomorphisms 
$\phi : R_1\to R_2$.

Let ${\hat\g}(t)\in (\bt\widehat\otimes\C[[t_{\bold H}]])^1$ be a 
solution to the 
Maurer-Cartan satisfying the condition of corollary 2.3. 
It follows from  this condition   
that the $\C[[t_{\bold H}]]$-module $T_{{\hat\g}(t)}$ is freely generated 
by the classes of partial 
derivatives $[\p_a{\hat\g}(t)]$. 
Therefore we have 
\proclaim{Proposition 3.1} There exists  formal power series  $A^c_{ab}(t)\in
\bt\widehat\otimes\C[[t_{\bold H}]]$ satisfying
$$ \p_a{\hat\g}\wedge\p_b{\hat\g}=\sum_c A^c_{ab}\p_c{\hat\g}\,\,\text{mod} 
\,\db_{{\hat\g}(t)}$$
The series  $A_{ab}^c(t)$ are the structure constants of the commutative 
associative \break $\C[[t_{\bold H}]]$@-algebra 
structure on  $\bold H\otimes\C[[t_{\bold H}]][-2]$.
\endproclaim
\qed
\Remark Note that on the tangent space at zero this algebra structure
is given by  the obvious multiplication on $\oplus_{p,q}H^q(M,\Lambda^pT_M)$.
This is "Mirror dual" to the ordinary multiplication on 
$\oplus_{p,q}H^q(M,\om^p_{\widetilde M})$.
\endremark

\section{4. Integral}
Introduce  linear functional on $\bt$ 
$$\int \g:=\int_M (\g\vdash\om)\wedge \om$$
\claim 4.1 It satisfies the following identyties:
 $$\aligned\int\db\g_1\wedge\g_2=(-1)^{\text{deg}\g_1}\int\g_1\wedge\db\g_2\\
\int\dv\g_1\wedge\g_2=(-1)^{\text{deg}\g_1+1}\int\g_1\wedge\dv\g_2
\endaligned\eqno (4.1)$$
for $\g_i\in  \Gamma(M,\Lambda^{q_i}\overline T^*_M\otimes\Lambda^{p_i}T_M), 
\,\, i=1,2$ where $\text{deg}\,\g_i=p_i+q_i-1$  .\endproclaim\qed

\section{5. Metric on $\Cal{T_M}$}

There exists a natural metric (i.e. a nondegenerate (super)symmetric 
$\Cal O_{M_\bt}$-linear pairing) on the sheaf $\Cal T_{\Cal M_\bt}$.
 In terms of a  solution to the Maurer-Cartan equation $\g^R
\in (\bt\otimes R)^1$ it 
means that there exists an $R$-linear graded symmetric pairing on $T_{\g^R}$,
which is functorial with respect to $R$. Here  $T_{\g^R}$
denotes the cohomology of $\db_{\g^R}$ defined in \S 3.  
The pairing is defined by the formula 
$$\langle h_1,h_2\rangle:=\int h_1\wedge h_2\,\,\,\text{for}\,\,h_1,h_2\in 
T_{\g^R} $$
where we assumed for simplicity that  $\g^R\in \Ker\dv\otimes R$.
It follows from (2.2) (see also lemma 6.1) that such a choice of $\g^R$ in the 
given class of gauge equivalence  is always possible. 
\claim 5.1 The pairing is compatible with the algebra 
structure.\endproclaim
\noindent{\qed}

\section{6. Flat coordinates on moduli space.} 

Another ingredient in the definition of Frobenius structure is the choice
of affine structure on the moduli space associated with $\bt$. 
The lemma 2.1 identifies $\Cal M_{\bt}$ with  the moduli space associated
with  trivial algebra $(\bold H[-1],d=0)$. The latter moduli space is 
the affine space $\bold H$. The affine coordinates $\{t_a\}$ on $\bold H$ 
give coordinates on $\Cal M_{\bt}$. This choice of coordinates on the 
moduli space 
corresponds to a specific
 choice of a universal solution to the Maurer-Cartan equation over 
 $\C[[t_{\bold H}]]$.
\proclaim{Lemma 6.1} There exists a  solution to the Maurer-Cartan equation in 
formal power series with values in $\bt$ 
$$\db{\hat\g}(t)+{[{\hat\g}(t),{\hat\g}(t)]\over 
2}=0,\,{\hat\g}(t)=\sum_a\hat\g_at^a+{1\over 
{2!}}\sum_{a_1,a_2}\hat\g_{a_1a_2}t^{a_1}t^{a_2} +\ldots\in
 (\bt\widehat\otimes\C[[t_{\bold H}]])^1$$ such that
\roster
\item (Universality) the cohomology classes $[\hat\g_a]$ form a basis of 
cohomology 
of the 
complex 
$(\bt,\db)$ 
\item (Flat coordinates) $\hat\g_a\in\Ker\dv,\,\,\hat\g_{a_1\ldots a_k}\in 
\im\dv\,\,\,
\text{for}\,\,\, k\geq 2$
\item (Flat identity) $\p_0\hat\g(t)=\bold 1$, where $\p_0$ is the coordinate 
vector field 
corresponding to $[\bold 1]\in \bold H[-1]$
\endroster
\endproclaim
\proof The theorem of \S 4.4 in [K2] shows that there exists $L_\infty$
 morphism $f$ 
homotopy inverse to the natural morphism $(\Ker\dv,\db)\to \bold H[-1]$ 
(for the definition of $L_\infty$-morphism see \S 4.3 in [K2]).  
Put $\Delta(t)=\sum_a(\Delta_a[-1]) t^a$ , where $\Delta_a[-1]$ denotes the
element $\dl_a$ having degree shifted by one. Then 
$${\g}(t)=\sum_n{1\over{n!}}f_n(\Delta(t)\wedge\dots\wedge\Delta(t))$$ 
satisfies 
the conditions $(1)-(2)$. To fulfill the condition $(3)$ $\g(t)$ must be 
improved slightly. Define the differential graded Lie algebra 
$\widetilde{\Ker}$
 as follows 
\roster
\item $\widetilde{\Ker}_i=\Ker \dv\subset \bt_i\,\text{for}\,i\geq 0$
\item $\widetilde{\Ker}_{-1}=\im \dv\subset \bt_{-1}$
\endroster
Note that the algebra $\Ker \dv$ is  the sum of the algebra 
$\widetilde{\Ker}$ and trivial algebra of constants $\C\otimes\bold 1[-1]$. 
Let 
$\tilde f$ be a homotopy inverse to the natural quasi-isomorphism 
$\widetilde{\Ker}\to \bold H[-1]_{\geq 0}$. Put $\tilde\Delta(t)=\sum_{a\neq 
0}(\Delta_a[-1]) t^a$. Then 
$$\hat{\g}(t)=\bold 1t_0 +\sum_n{1\over{n!}}\tilde 
f_n(\tilde\Delta(t)\wedge\dots\wedge\tilde\Delta(t))$$ satisfies all
the conditions. 
\endproof
\Remark Any two formal power series  satisfying conditions of lemma 6.1 are 
equivalent under the natural action of the group associated with the \break Lie 
algebra  
$(\widetilde\Ker \dv\widehat\otimes\C[[t_{\bold H}]])^0$. 
\endremark
\Remark It is possible to write down an explicit formula for the components of 
the morphism $f$ in terms of Green functions of the Laplace operator acting 
on differential forms on $M$.
\endremark
\Remark  After the identification  of the moduli space $\Cal M_\bt$ with 
$\bold H$, provided by lemma 2.1, the complex structure moduli space 
corresponds to the \break subspace  $\bold H^1(M,\Lambda^1T_M)$.
In the case of classical moduli space of the complex structures on $M$
the analogous lemma was proved in [T]. The coordinates arising on 
 the classical moduli 
space of complex structures coincide with so called 
"special" coordinates of [BCOV].\endremark
\remark{Notation}
Denote $${\hat\g}(t)=\sum_a\hat\g_at^a+{1\over 
{2!}}\sum_{a_1,a_2}\hat\g_{a_1a_2}t^{a_1}t^{a_2} +\ldots\,\in (\bt\otimes\bold 
C[[t_{\bold H}]])^1$$
a solution to the Maurer-Cartan equation satisfying conditions of lemma 6.1.
\endremark

 The parameters of a miniversal solution to the Maurer-Cartan equation 
over $\C[[t_{\bold H}]]$ serve as  coordinates on the moduli space. The 
specific choice of coordinates corresponding to   the solution to the 
Maurer-Cartan equation satisfying conditions (1)-(2) of lemma 6.1  corresponds 
to choice
of  coordinates  on moduli space that are flat with
respect to the natural (holomorphic) metric $g_{ab}$.

\claim 6.2  The power series $\langle 
\p_a{\hat\g}(t),\p_b{\hat\g}(t)\rangle\in\C[[t_{\bold H}]]$ has only  
constant term in the  power series expansion at $t=0$.\endproclaim 
\proof $\langle x,y\rangle=0$
for $x\in \Ker\dv, y\in \im\dv$.\endproof

\remark{Notation}
Denote  $g_{ab}:=\langle \p_a{\hat\g}(t),\p_b{\hat\g}(t)\rangle$.
\endremark
Thus we have constructed all the ingredients of the Frobenius structure on $\Cal 
M_{\bt}$:
the tensors $A^c_{ab}(t)$, $g_{ab}$ and the coordinates $\{t_a\}$ that are flat 
with
 respect
to $g_{ab}$.
\Remark The 3-tensor $A_{ab}^c(t)$ on $\Cal M_\bt$ does not depend on the choice
of $\om$. The 2-tensor $g_{ab}$ is multiplied by $\lambda^2$ when $\om$ is 
replaced by
$\lambda\om$\endremark

\claim 6.3 The structure constants satisfy $A_{0a}^b=\delta^b_a$.\endproclaim
\proof It follows from the condition (3) imposed on $\hat\g(t)$\endproof

We have checked that the tensors $A^c_{ab}(t),\,g_{ab}$ have the properties 
(1)-(3)  
from the definition of the Frobenius structure. It remains to us to check the
\break property (4).  

\section{ 7. Flat connection and periods}

Let ${\hat\g}(t)\in (\bt\otimes\bold 
C[[t_{\bold H}]])^1$ be
a solution to the Maurer-Cartan equation satisfying conditions (1)-(2) of lemma 
6.1.
Then the formula
$$\om(t):=e^{\hat\g(t)}\vdash\om $$ defines a closed form of mixed degree
depending formally on $t\in \bold H$. 
For $t\in \bold H^{-1,1}$ $\hat\g(t)\in \Gamma(M,\overline T\otimes T^*)$ 
represents a 
deformation of complex 
structure.  Then $\om(t)$ is
a  holomorphic $n$-form in the complex structure corresponding to $t\in \bold 
H^{-1,1}$, 
where $n=\text{dim}_\C\, M$.

Let $\{p^a\}$ denote the set of sections of  
$\Cal T^*_{\bold H}$ that form a framing dual to $\{\p_a\}$.
Define a (formal) connection on $\Cal T^*_{\bold H}$ by the covariant 
derivatives:
$$\nabla_{\p_a}(p^c):=\sum_b A_{ab}^c p^b \eqno (7.1)$$
Strictly speaking this covariant derivatives are formal power series
sections of $\Cal T^*_{\bold H}$.

Let us put 
$$\Pi_{ai}={\p\over\p t^a}\int_{\Gamma_i} \om(t) \eqno 
(7.2)$$
where $\{\Gamma_i\}$ form a basis in $\bold H_*(M,\C)$. In particular
$$\Pi_{0i}=\int_{\Gamma_i} \om(t)$$
if $\hat\g(t)$ satisfies the condition (3) of lemma 6.1.
\lemma 7.1 The periods
$\Pi_{i}=\sum_{a}\Pi_{ai} p^a$ are flat sections of $\nabla$\endproclaim
\noindent{\it Proof.\,}Let $\p_\ta=\sum_a \ta^a \p_a$ be an even constant vector 
field, i.e.  $\ta^a$ are even constants for even $\p_a$ and odd for odd $\p_a$. 
It is enough to prove that 
$$\p_\ta\p_\ta\int_{\Gamma_i}\om(t)=\sum_c A^c_{\ta\ta}\p_c \int_{\Gamma_i}
\om(t)$$
where $A^c_{\ta\ta}$ are the algebra structure constants defined via 
$\sum_a\ta^a\p_a
\circ\sum_a\ta^a\p_a=\sum_c 
A^c_{\ta\ta}\p_c$ (see \S 1).
Note that the operators $\dv$ and $\db_{\hat\g(t)}$ acting on 
$\bt\widehat\otimes\C[[t_{\bold H}]]$
satisfy a version of $\p\db$-lemma :
$$\Ker 
\dv\cap\Ker\db_{\hat\g(t)}\cap(\im\dv\oplus\im\db_{\hat\g(t)})=
\im\dv\circ\db_{\hat\g(t)}\eqno(7.3)$$
Equivalently, there exists decomposition of 
$\bt\widehat\otimes\C[[t_{\bold H}]]=X_0\oplus X_1\oplus
X_2\oplus X_3\oplus Y$ into direct sum of 
 graded vector spaces, such that
  the only nonzero components of $\db_{\hat\g(t)}$, $\dv$ are isomorphisms
   $\db_{\hat\g(t)}:X_0\mapsto X_1,\,\,\,X_2\mapsto X_3;\,\,\,\,\,
 \dv:X_0\mapsto X_2,\,\,\,X_1\mapsto X_3$.
 Differentiating twice the Maurer-Cartan equation (2.4)
with respect to $\p_\ta$ and using (2.1) one obtaines
$$\dv (\p_{\ta}\hat\g(t)\wedge\p_{\ta}\hat\g(t)-\sum_c 
A_{\ta\ta}^c\p_c\hat\g(t))=-\db_{\hat\g(t)}\p_{\ta}\p_{\ta}\hat\g(t) \eqno 
(7.4)$$
It follows from $\p\db$-lemma for $\dv$, $\db_{\hat\g(t)}$  and the equation 
(7.4) 
that  there exist formal power series  
$\alpha_{\ta}(t)\in \bt\widehat\otimes 
\C[[t_{\bold H}]]$ such that 
$$\eqalign{&\p_{\ta}\hat\g(t)\wedge\p_{\ta}\hat\g(t)-\sum_c 
A_{\ta\ta}^c\p_c\hat\g(t)=\db_{\hat\g(t)}\alpha_{\ta}(t),\cr  & 
\p_{\ta}\p_{\ta}\hat\g(t)=\dv(\alpha_{\ta}(t)) \cr}$$ 
Therefore
$$\p_{\ta}\p_{\ta}\om(t)=\sum_c 
A_{\ta\ta}^c{\p_c}\om(t) + 
d(\alpha_{\ta}e^{\hat\g(t)}\vdash\om)$$\qed

 It follows from the condition (1) imposed on $\hat\g(t)$ that $\Pi_i$ form
 a (formal) framing of $\Cal T^*_{\bold H}$.
 \corollary 7.2 The connection $\nabla$ is  flat\endproclaim

\claim 7.3 The structure constants $A^c_{ab}$ satisfy the potentiality condition 
(4)
 in flat coordinates.\endproclaim
\proof
If one puts  symbolically $\nabla=\nabla_0 + A$ then the flatness of $\nabla$
implies that
$$\nabla_0 A+{1\over 2}[A,A]=0$$ 
Notice that
associativity and commutativity of the algebra defined by $A^c_{ab}$ imply
that $$[A,A]=0 .$$
Therefore
$$\nabla_0(A)=0. $$
\endproof

We have completed the proof of the fact that $A^c_{ab}(t)$ and $g_{ab}$
 define the Frobenius structure  on $\Cal M_\bt$ in the flat coordinates 
$\{t_a\}$.
  
\Remark In fact one can write an explicit formula for the potential of the 
Frobenius structure. Let us put ${\hat\g}(t)=\sum_a\hat\g_a t^a+\dv\alpha(t), 
\alpha(t)\in (\bt\widehat\otimes  t_{\bold H}^2\C[[t_{\bold H}]])^0$. Put 
$$\Phi=\int 
 -{1\over 2}\db\alpha\wedge\dv\alpha+{1\over 
6}{\hat\g}\wedge{\hat\g}\wedge{\hat\g}$$
Then one checks easily that $A_{abc}=\p_a\p_b\p_c\Phi$ (see Appendix). 
In the case $\text{dim}_\C M=3,\, \hat \g \in 
\bt^{-1,1}=\Gamma(M,\overline{T}^*_M\otimes T_M)$
this formula gives the critical value of so called Kodaira-Spencer Lagrangian of 
[BCOV].
\endremark

\Remark Define differential Batalin-Vilkovisky algebra as  $\Z_2$-graded
commutative associative
algebra $A$ equipped with odd differentiation $\db,\,\, \db^2=0$ and odd 
differential operator $\dv$ of order $\leq 2$ such that
$\dv^2=0,\,\,\,\dv\db+\db\dv=0,\,\,\,\dv(1)=0$. One can use the formula (2.1)  
to define  the structure of $\Z_2$-graded
Lie algebra on $\Pi A$.  Assume that the operators $\db,\dv$ satisfy 
$\p\db$-lemma (2.3).
Assume in addition that $A$ is equipped with a linear functional  
$\int:A\to \C$ 
satisfying (4.1) such that the metric defined as in \S 5 is nodegenerate.  Then 
the same construction as above produces the Frobenius  
structure 
on the $\Z_2$-graded moduli space $\Cal M_{\Pi A}$.  One can define
the tensor product of two such Batalin-Vilkovisky algebras.  
Operator $\dv$ on $A_1\otimes A_2$ is given by 
$\dv_1\otimes 1+1\otimes\dv_2$. Also, $\db$  on  $A_1\otimes A_2$ is 
$\db_1\otimes 1+1\otimes\db_2$.
It is naturally to expect that the Frobenius manifold corresponding to 
$A_1\otimes A_2$
is equal to the tensor product of Frobenius manifolds corresponding to 
$A_1,A_2$, defined in [KM] in terms of the corresponding algebras over operad 
$\{H_*(\overline
M_{0,n+1})\}$.
\endremark

\section{8. Scaling transformations}

The vector field $E=\sum_{a}-{1\over 2}|\dl_a|t^a\p_a$ generates the
scaling symmetry on $\bold H$. 
\proclaim{Proposition 8.1} ${\Cal Lie}_EA_{abc}=({1\over 
2}(|\dl_a|+|\dl_b|+|\dl_c|)+3-\text{dim}_\C M)A_{abc}$ 
\endproclaim
\proof  $A_{abc}=\int_M \p_{a}{\hat\g}\wedge\p_{b}{\hat\g}\wedge\p_{c}{\hat\g}$.  
Note that
 $\int \g\neq 0$ implies that $\g\in \bt_{2n-1}$. The proposition follows from
the grading condition on $\hat\g(t)$.
\endproof

\corollary 8.2 ${\Cal Lie}_EA_{ab}^c=({1\over 
2}(|\dl_a|+|\dl_b|-|\dl_c|)+1)A_{ab}^c$ 
\endproclaim 

Note that the vector field $E$ is conformal with respect to the metric 
$g_{ab}$.
Therefore
the proposition 8.~1 shows that $E$ is the Euler vector field of the  Frobenius 
structure on $\bold H$ (see [M]). Such a vector field is defined uniquely up to
a multiplication by a constant. 
The simplest invariant of
 Frobenius manifolds is the spectrum
of the operator $[E,\cdot]$ acting on infinitesimal generators of translations
 and the weight of the
tensor $A_{ab}^c$ under the ${\Cal Lie}_E$-action.
Usually  
the normalization of $E$ is chosen so that $[E,\p_0]=1$.

 In our case the spectrum of $[E,\cdot]$ is equal to   $$\bigcup_d \{1-d/2\} 
\,\,\text{\eightrm with  multiplicity}\sum_{q-p=d-n}\text{dim}H^q(M,\Omega^p)$$
Note that this spectrum and the weight of   $A_{ab}^c$ coincide identically 
with the corresponding quantities of the Frobenius structure arising from the 
Gromov-Witten invariants of the dual Calabi-Yau manifold  $\widetilde M$.

\section {9.Further developments.}
Conjecturally the constructed Frobenius manifold is related to the Gromov-Witten 
invariants of
$\widetilde M$
in the following way. One can rephrase the present construction 
in purely algebraic terms using \v Cech instead of Dolbeault realization of the 
simplicial 
graded Lie algebra $\Lambda^*\Cal T_M$. The only
additional ``antiholomorphic'' ingredient that is used is the choice of a 
filtration
on $H^*(M,\C)$  complementary to the Hodge filtration. The Frobenius structure, 
arising from the limiting weight filtration corresponding to a point with 
maximal 
unipotent
monodromy  on moduli space of complex structures on M, 
coincides conjecturally  with
 Frobenius structure on $H^*(\widetilde M,\C)$, obtained from the
  Gromov-Witten invariants.

Our construction of Frobenius manifold  is a particular case of 
a more general construction. Other cases of this construction include
the Frobenius manifold structure on the moduli space of singularities of 
analytic
functions found by K.\,Saito, the Frobenius manifold structure on the
moduli space of ``exponents of algebraic functions''. The latter case is a 
Mirror Symmetry
partner to the structure arising from Gromov-Witten invariants on Fano 
varieties.
All these cases seem to be related with yet undiscovered generalization of 
theory of 
Hodge structures.
We hope to return to this in the next paper.

\specialhead Appendix \endspecialhead

Let ${\hat\g}(t)=\sum_a\hat\g_a t^a+\dv\alpha(t), \alpha(t)\in 
({\bold t}\otimes  t_{\bold H}^2{\bold C}[[t_{\bold H}]])^1$ be
a solution to Maurer-Cartan equation satisfying conditions (1)-(2) of lemma 
6\.1.
Put $$\Phi=\int  -{1\over 2}\db\alpha\wedge\dv\alpha+{1\over 
6}{\hat\g}\wedge{\hat\g}\wedge{\hat\g}$$ 
\proclaim{Proposition } $A_{abc}=\p_a\p_b\p_c\Phi$
\endproclaim
\proof. Let $\p_{\ta}=\sum_a \ta_a\p_a$ be an even constant vector 
field in $\bold H$. It is enough to prove
that $\int 
(\p_{\ta}{\hat\g}\wedge\p_{\ta}{\hat\g}\wedge{\p_{\ta}\hat\g})=
\p^3_{\ta\ta\ta}\Phi$.  Let us differentiate the terms in $\Phi$ 
$$\eqalign{&\p^3_{\ta\ta\ta}({\hat\g}\wedge\hat\g\wedge\hat\g)=
18\p^2_{\ta\ta}{\hat\g}
\wedge\p_{\ta}{\hat\g}\wedge{\hat\g}+
3\p^3_{\ta\ta\ta}{\hat\g}\wedge{\hat\g}\wedge\hat\g+
6\p_{\ta}{\hat\g}\wedge\p_{\ta}{\hat\g}
\wedge\p_{\ta}{\hat\g} \cr}$$
$$\eqalign{&\p^3_{\ta\ta\ta}(\db\alpha\wedge\dv\alpha)=
\db\alpha\wedge(\p^3_{\ta\ta\ta}\dv\alpha)+
(\p^3_{\ta\ta\ta}\db\alpha)\wedge\dv\alpha 
+3(\p_{\ta}\db\alpha)\wedge(\p^2_{\ta\ta}\dv\alpha)+\cr
&+3(\p^2_{\ta\ta}\db\alpha)\wedge(\p_{\ta}\dv\alpha)\cr}\eqno(*)$$

Notice that 
$$\eqalign{&\int(\p^3_{\ta\ta\ta}\db\alpha)\wedge\dv\alpha=
(-1)^{\text{deg}\db\alpha+1}\int(\p^3_{\ta\ta\ta}\dv\db\alpha)\wedge\alpha=
\int(\p^3_{\ta\ta\ta}\dv\db\alpha)\wedge\alpha=\cr
&=-\int(\p^3_{\ta\ta\ta}\db\dv\alpha)\wedge\alpha=
-(-1)^{\text{deg}\dv\alpha}\int(\p^3_{\ta\ta\ta}\dv\alpha)\wedge\db\alpha=
\int(\p^3_{\ta\ta\ta}\dv\alpha)\wedge\db\alpha=\cr
&=(-1)^{(\text{deg}\dv\alpha+1)(\text{deg}\db\alpha+1)}\int
\db\alpha\wedge(\p^3_{\ta\ta\ta}\dv\alpha)=
\int\db\alpha\wedge(\p^3_{\ta\ta\ta}\dv\alpha)\cr}$$

\noindent Hence, the first two terms in (*) give the same contribution.

\noindent We have $$\eqalign{&\int\db\alpha\wedge(\p^3_{\ta\ta\ta}\dv\alpha)=
(-1)^{\text{deg}\db\alpha+1}\int\dv\db\alpha\wedge(\p^3_{\ta\ta\ta}\alpha)
=\int\dv\db\alpha\wedge(\p^3_{\ta\ta\ta}\alpha)=\cr
&=-\int\db\dv\alpha\wedge(\p^3_{\ta\ta\ta}\alpha)=
-\int\db\hat\g\wedge(\p^3_{\ta\ta\ta}\alpha)=
{1\over 2}\int[\hat\g,\hat\g]\wedge(\p^3_{\ta\ta\ta}\alpha)=\cr
&={1\over 2}\int\dv(\hat\g\wedge\hat\g)\wedge(\p^3_{\ta\ta\ta}\alpha)=
(-1)^{(\text{deg}(\hat\g\wedge\hat\g)+1)}{1\over 
2}\int(\hat\g\wedge\hat\g)\wedge(\p^3_{\ta\ta\ta}\dv\alpha)=\cr
&={1\over 2}\int(\hat\g\wedge\hat\g)\wedge\p^3_{\ta\ta\ta}\hat\g\cr}$$

\noindent Similarly,$$\eqalign{ 
\int(\p_{\ta}\db\alpha)\wedge(\p^2_{\ta\ta}\dv\alpha)&={1\over 2}\int
\p_{\ta}(\hat\g\wedge\hat\g)\wedge\p^2_{\ta\ta}\hat\g \cr
\int(\p^2_{\ta\ta}\db\alpha)\wedge(\p_{\ta}\dv\alpha)&={1\over 2}\int
(\p^2_{\ta\ta}\hat\g)\wedge\p_{\ta}(\hat\g\wedge\hat\g) \cr
}$$
\endproof

\enddocument
\end